\documentclass[twocolumn]{svjour3}         

\smartqed  

\usepackage{graphicx}
\usepackage{mathptmx}
\usepackage{latexsym}

\usepackage{algorithmic}
\usepackage{algorithm}
\usepackage{amsmath,amssymb}
\usepackage{multirow}
\usepackage{subfigure}

\usepackage{makeidx}         
\usepackage{multicol}        
\usepackage[bottom]{footmisc}

\begin{document}

\renewcommand\b{{\bf b}}
\newcommand\f{{\bf f}}
\renewcommand\t{{\bf t}}
\newcommand\xx{{\bf x}}
\newcommand\XX{{\bf X}}

\newcommand\bz{{\bf 0}}
\newcommand\commentout[1]{\relax}
\newcommand\R{{\bf R}}
\renewcommand\u{{\bf u}}
\newcommand\X{{\bf X}}
\newcommand\x{{\bf x}}
\newcommand\eref[1]{$(\ref{#1})$}
\newcommand\s{{\bf s}}
\renewcommand{\thefootnote}{\fnsymbol{footnote}} 
\newcommand{\bo}{\mathbf}
\renewcommand{\b}{\bo{b}}
\renewcommand{\c}{\bo{c}}
\renewcommand{\d}{\bo{d}}
\newcommand{\e}{\bo{e}}
\newcommand{\g}{\bo{g}}
\newcommand{\h}{\bo{h}}
\renewcommand{\i}{\bo{i}}
\renewcommand{\r}{{\bo r}}
\newcommand\atsign{{@@}}
\let\vaccent\v
\renewcommand{\v}{{\bo v}}
\newcommand\trace{\mbox{trace}}
\newcommand\Grad{\mbox{Grad}}
\newcommand\area{\mbox{area}}
\newcommand\volume{\mbox{volume}}
\newcommand\diag{\mbox{diag}}
\newcommand\tril{\mbox{tril}}
\newcommand\triu{\mbox{triu}}
\newcommand{\w}{{\mathbf w}}
\newcommand{\y}{{\mathbf y}}
\newcommand{\z}{{\mathbf z}}
\newcommand\chiab{\bar\chi_A}
\newcommand\chia{\chi_A}
\newcommand\chiamub{\bar\chi_{A,\mu}}
\newcommand\chiamu{\chi_{A,\mu}}
\newcommand\eps{\epsilon}
\newcommand\C{{\mathbf C}}
\newcommand\sref[1]{Section \ref{#1}}
\newcommand\thref[1]{Theorem \ref{#1}}
\newcommand\lemref[1]{Lemma \ref{#1}}
\newcommand\tabref[1]{Table \ref{#1}}
\newcommand\cD{{\cal D}}
\newcommand\Mone{\ensuremath{A_{1}^{T}D_1 A_1}}
\newcommand\mone{\ensuremath{A_{1}^{T}D_1}}
\newcommand\Mtwo{\ensuremath{A_{2}^{T}D_2 A_2}}
\newcommand\mtwo{\ensuremath{A_{2}^{T}D_2}}
\renewcommand\Re{\mathop{\rm Re}}
\renewcommand\Im{\mathop{\rm Im}}

\title{A robust solution procedure for hyperelastic solids with large boundary 
deformation\footnote{Much of the work of the first author
was performed while a graduate student in the Center for
Applied Mathematics at Cornell University.  The author's
work is supported by a fellowship from the National Physical 
Science Consortium (with support from Sandia National Laboratories, and 
Cornell University); it is also supported in part by NSF grants 
ACI-0085969, CNS-0720749, and NSF CAREER Award OCI-1054459.}\footnote{
The work of the second author was supported in part by NSF grant 
ACI-0085969, a Discovery grant from 
NSERC (Canada), and a grant from the U.S. Air Force Office of 
Scientific Research.}\footnote[3]{Paper accepted for publication in 
{\emph{Engineering with Computers}}.  The final publication is available 
at www.springerlink.com (DOI:  10.1007/s00366-011-0225-y).}}  

\titlerunning{Solution Procedure for Hyperelastic Solids with Large Deformation}        

\author{Suzanne M. Shontz \and
        Stephen A. Vavasis 
}


\institute{Suzanne M. Shontz\at
                Department of Computer Science and Engineering \\
		The Pennsylvania State University\\
              University Park, PA 16802, United States\\
              \email{shontz@cse.psu.edu}           
           \and
           Stephen A. Vavasis \at
           Department of Combinatorics and Optimization \\
           University of Waterloo\\
           Waterloo, Ontario, N2L 3G1, Canada \\
		\email{vavasis@math.uwaterloo.ca}
}

\date{Received: date / Accepted: date}

\maketitle

\begin{abstract}
Compressible Mooney-Rivlin theory has been us-\\ed
to model hyperelastic solids, such as rubber and porous polymers,
and more recently for the modeling of soft tissues for biomedical 
tissues, undergoing large elastic deformations.
We propose a solution procedure for Lagrangian finite element
discretization of a static nonlinear compressible Mooney-Rivlin
hyperelastic solid.
We consider the case in which the boundary
condition is a large prescribed deformation, so that mesh
tangling becomes an obstacle for straightforward
algorithms.  Our solution procedure involves a largely
geometric procedure to untangle the mesh: solution of a sequence of linear
systems to obtain initial guesses for interior nodal
positions for which no element is inverted.  After the
mesh is untangled, we take Newton iterations to converge to a
mechanical equilibrium.  The Newton iterations are safeguarded
by a line search similar to one used in optimization. Our
computational results indicate that the algorithm is up
to 70 times faster than a straightforward Newton 
continuation procedure and is also more robust (i.e., able to tolerate 
much larger deformations).  For a few extremely large deformations, the
deformed mesh could only be computed through the use of an 
expensive Newton continuation method while using a tight convergence 
tolerance and taking very small steps. 
\end{abstract}

\section{The problem under consideration}
\label{sec:intro}

We consider the problem of
solving for the deformed shape of a hyperelastic solid
body under static loads.
The continuum mechanical model under consideration has the following
description \cite{Holzapfel}.  Let $B_0\subset \R^d$ be an undeformed 
solid body
whose boundary is $\partial B_0$.  Here $d$, the space dimension, is 2 
or 3. Assume boundary conditions (either displacement or traction, i.e.,
Dirichlet or Neumann) are given as follows.  The boundary $\partial B_0$
is partitioned into two subsets $\Gamma_D$ and $\Gamma_N$.  A function
$\phi_0\mbox{:}\Gamma_D\rightarrow \R^d$ specifies new-position 
(Dirichlet)
boundary conditions.  A second function $\t_0\mbox{:}\Gamma_N\rightarrow 
\R^d$ 
specifies
traction (Neumann) boundary conditions. 
Everything in this paper 
extends to the more general case that some coordinate entries are
Neumann while others are Dirichlet at certain boundary points, but we
limit the
discussion to the special case that each boundary point
is Dirichlet or Neumann in all $d$ coordinates
in order to simplify notation.
Finally, the model requires a 
specification of the model's body forces, 
that is, a function $\b\mbox{:}B_0\rightarrow \R^d$ that specifies
the force of gravity and other forces on the body.

The problem is to find a function $\phi\mbox{:}B_0\rightarrow {\R}^d$ 
that
specifies the new position of the body.  Let $B$ denote $\phi(B_0)$.
For a point $\XX\in B_0$, let $\x=\phi(\XX)$.  Let $F$ be the deformation 
gradient, i.e., $F=d\phi/d\XX=d\x/d\XX$.  
It is assumed that $F(\XX)$ has a positive 
determinant for all $\XX$. 
The {\em Green-Lagrange strain tensor} is defined to be $E=(F^TF-I)/2$. 
Let scalar function $\Psi(F)$ be the {\em strain energy function}, which 
is assumed to be a property of the material. For this paper, we assume 
that $\Psi$ depends only on
two scalar invariants of tensor $E$,
namely 
$J=\det(F)=\sqrt{\det(2E+I)}$ and
$I_1=\trace(F^TF)=\trace(2E+I)$.
Further specializing the model, the strain energy is then 
taken to have the following 
form suggested by Ciarlet and Geymonat in \cite{Ciarlet} for compressible 
Mooney-Rivlin 
materials \begin{equation} 
\Psi(F) = 
\frac{\lambda}{4}(J^2-1)-\left(\frac{\lambda}{2}+\mu\right)\ln 
J+\frac{\mu}{2}(I_1-3)\label{eq:MR-model} \end{equation} where 
$\lambda, \mu >0$ are material parameters.
Compressible Moo-\\ney-Rivlin theory has been used for analyzing large 
elastic deformations of soft materials, including 
rubber~\cite{mooney_rivlin_original}; porous polymers, such as
porous polyethylenes used as insulation boards for construction, 
protective packaging materials, insulated drinking cups, and flotation 
devices~\cite{gibson_ashby}; and biological tissues~\cite{weichert};
as well as other applications.

For the $d = 2$ case, we assume that $B_0$ is 3D but that the 
$z$-displacement is identically $0$ and that the $x$- and $y$-displacements 
depend only on $x$ and $y$; these are called
{\em  plane strain} assumptions.  
Thus, $E$ has a last row and column of all zeros, and the Mooney-Rivlin
formula in \eref{eq:MR-model} is applied to this $E$ to come up with 
the strain energy function for the $d = 2$ case.   
The condition for static equilibrium (written in minimization form) 
is that
\begin{equation} \int_{B_0} 
\Psi(F(\XX))\,dV-\int_{B_0}\rho\b\cdot\phi(\XX)\,dV
-\int_{\Gamma_N}\t_0\cdot\phi(\XX)\, dA
 \label{eq:minform}
\end{equation} 
is minimized among all choices of $\phi$ that satisfy the
Dirichlet boundary condition, i.e., that satisfy $\phi(\XX)=\phi_0(\XX)$
for all $\XX\in \Gamma_D$.

This
condition can be rewritten in variational form: for all admissible 
variations $\delta\u$, that is, functions in the space $[H^1(B_0)]^d$ 
that vanish on $\Gamma_D$, \begin{equation} \int_{B_0} \frac{\partial 
\Psi}{\partial F} \mbox{:} \Grad \ \delta \u \ dV - \int_{B_0}
\rho\b\cdot\delta\u\,dV-\int_{\Gamma_N}\t_0\cdot\delta\u\, dA
= 0,
\label{eq:varform} \end{equation}
where $A\mbox{:}B = \trace(AB^T)$ is used to denote the inner product of 
second-order tensors $A$ and $B$.  
This model also applies to the case of linear 
elasticity with two changes in definitions.  
First, $E = ((F-I)^T+(F-I))/2$ in the
case of linear elasticity.  Second, 
$\Psi = \mu \sum_{i,j} E(i,j)^2+\frac{\lambda}{2} 
\left( \sum_{i} E(i,i) \right)^2$, which can be written in terms of 
the two 
invariants of $E$.  

We should mention that our method does not appear to depend so 
much on the specific details of the Mooney-Rivlin model, except for the
$\ln J$ term, which is quite important for our analysis.  Since 
$dv=J\,dV$, where $dv$ is the volume element of $B$ and $dV$ is the
volume element of $B_0$, this logarithmic term resists infinite 
compression of the material: if a small positive volume of material in 
$B_0$ shrinks to a 0-volume set in $B$, then this term causes the strain 
energy at those points to become infinite.

We next describe the Lagrangian 
discretization of the problem under consideration
\cite{BelytschkoLiu}.  We assume that $B_0$ is discretized with a mesh of
triangles or tetrahedra.  We assume that the discretization of $\phi$, or
alternatively the discretization of the displacement $\u=\phi(\XX)-\XX$, is
piecewise linear, with the pieces of linearity being the mesh cells. (In
Section~\ref{sec:conclusions}, we discuss extension of our method to
piecewise quadratic displacements.)  Recall that $d$ is the space 
dimension, and let $m$ denote the number of non-Dirichlet nodes of the 
mesh.  This assumption implies that $\u$ is
determined by $dm$ real numbers, namely, the values of $\u$ at nodes.
The finite element method finds the
displacement $\u$ such that \eref{eq:varform} holds for all test functions 
$\delta\u$ in the test function space.  Here, the 
test function space is the set of $\delta\u$'s that are piecewise linear 
and
continuous and vanish on $\Gamma_D$. The integral in
\eref{eq:varform} is evaluated with a quadrature rule; we have used a
6-point formula having degree 4 precision from \cite{Hughes} for our
quadrature in 2D and a 15-point formula having degree 5 precision from
\cite{Keast} for 3D.  
It suffices to
solve \eref{eq:varform} for the $dm$ choices of $\delta\u$ that compose
the standard basis for the test function space.  This yields a system of
$dm$ nonlinear equations for $dm$ unknowns.

The algorithmic question under consideration is how to robustly solve
these nonlinear equations.  In the next section, we give a summary of
the mesh tangling issue and of our proposal to overcome it.
The
individual steps of our algorithm are then described in more detail in
Sections~\ref{sec:femwarp} and
\ref{sec:NewtonLineSearch}.  
In Section \ref{sec:contin}, we summarize the Newton continuation 
algorithm which is a popular technique within the engineering 
community for solving the nonlinear equations.
Our computational experiments, which compare the two algorithms, are 
presented in Sections~\ref{sec:2Dexperiments} and \ref{sec:3Dexperiments}.
Concluding remarks, including
some discussion of the incompressible case,
are presented in Section~\ref{sec:conclusions}.

The preceding formulation is called ``Lagrangian'' discretization
because the nodes of the mesh remain fixed with respect to material
points throughout the solution procedure.
Alternatives to the Lagrangian approach
include the Eulerian approach
and arbitrary Lagrangian-Eulerian (ALE) methods.  Pure Eulerian methods
are not widely used in solid mechanics because of the difficulty in
applying boundary conditions.  ALE methods are a more viable competitor
to Lagrangian methods; in ALE methods the geometry is reme-\\shed as part
of the solution procedure.  ALE remeshing attempts to preserve a high-quality
mesh as the solution evolv-\\es.  ALE methods are substantially more complicated
than Lagrangian methods because of the need to interpolate field quantities
to new mesh points on every remeshing step.  In addition, ALE remeshing is
itself somewhat of an art in that there is no foolproof universal procedure for
updating the mesh.

For these reasons, we focus on traditional Lagrangian solution
techniques in this paper.  Nonetheless,
the first part of our algorithm (called ``iterative stiffening'' in
Section~\ref{sec:femwarp}) can be regarded as a particular ALE remeshing
approach; we return to this topic later.

\section{Mesh tangling}
\label{sec:tangling}

The standard method for solving a system of nonlinear equations is
Newton iteration.  It is well-known, however, that if the initial
guess is far from the true solution, then Newton iteration will
often diverge.  

In the case of hyperelasticity with large deformation, there is a specific
obstacle that may cause divergence, namely, mesh tangling.  The definition 
of this term is that a mesh is {\em tangled} if the value of $J$ defined 
in the previous section is 0 or negative in $B_0$. In the case of linear 
displacements, $J$ is piecewise constant, and hence this condition can 
be 
verified with a finite number of determinant computations.  The matter of 
checking for tangling in the piecewise quadratic case is more complicated 
and is discussed in Section~\ref{sec:conclusions}.  A solution with a 
tangled mesh is physically invalid.  Indeed, the strain energy function is 
undefined in this case because of the presence of the term 
$-\left(\frac{\lambda}{2}+\mu\right)\ln J$.
Note that although the strain energy function is undefined
when $J$ is negative, the Galerkin form \eref{eq:varform} is
still well defined, which is an anomaly that we return to below.
We assume that the given problem instance
has a valid solution, i.e., there is a piecewise linear function
$\u$ satisfying the boundary conditions, as well as
\eref{eq:varform} for all
test functions $\delta\u$ plus  the condition that $J>0$ on every element.

Even with this assumption, Newton's method will still often run into
problems because the mesh will become tangled on intermediate steps.
For example, the starting point for Newton's method is often taken to
be $\u=\bz$ on every interior node.  If the deformation of the boundary is
large, then this starting point corresponds to a mesh which will have 
tangling among most of the elements that are adjacent to the boundary.

To understand a difficulty posed by a tangled mesh, suppose that the 
strain energy has a single term 
$$\Psi(J)=-\left(\frac{\lambda}{2}+\mu\right)\ln J$$
(for $\frac{\lambda}{2}+\mu > 0)$ 
on a single element and that there are no boundary constraints.  If we 
treat $J$ as the independent scalar variable, then Newton's method for 
minimizing this scalar function is 
$$J^{(i+1)}=J^{(i)}-\Psi'(J^{(i)})/\Psi''(J^{(i)})$$ which simplifies in 
this case to $$J^{(i+1)}=2J^{(i)}.$$ For positive $J^{(0)}$, this 
iteration produces a sequence of $J$'s tending to $+\infty$. This is to be 
expected since the minimum of $-\left(\frac{\lambda}{2}+\mu\right)\ln J$ is
indeed
at $+\infty$. On the other hand, for a negative $J^{(0)}$, this iteration
tends to $-\infty$, which is physically invalid.

The preceding analysis, although naive, seems to point
to the following conclusion:
Newton's method on the Galer-\\kin form, when applied to a tangled
mesh, has a natural tendency to make the tangling worse.  
We suspect that this fact is probably already known to experts in the
field, although we have not been able to find it in the previous literature.

Given the conclusion in the previous paragraph, it seems of paramount
importance to avoid tangling.  When Newton's method fails in computational
mechanics, it is standard practice to try Newton continuation,
that is, to apply the load in incremental steps and use the converged 
solution for one step as the Newton starting point for the next step. 
Continuation is described in more detail in Section~\ref{sec:contin}.
Continuation, however, addresses the 
tangling issue only in an indirect fashion and therefore is likely to be 
very inefficient.  Our computational experiments confirm the inefficiency 
of continuation.

We propose a new algorithm for getting around the mesh tangling obstacle.
The basic idea is to first untangle the mesh using a much simplified
mechanical model.
Once the mesh is untangled, the true mechanical model is solved.
``Untangling the mesh'' means
finding a
$\phi$ that satisfies the Dirichlet boundary condition and also
satisfies $J>0$.
Our new algorithm, which we call UBN (for ``untangling before Newton'') 
consists of two steps.
\begin{enumerate}
\item
First, we attempt to untangle the mesh with the 
iterative-stiffening algorithm, described in Section~\ref{sec:femwarp}.
Iterative stiffening builds on the
FEMWARP algorithm from
our previous work \cite{ShontzVava}.  That paper, however, concerned 
itself with a pure mesh generation problem
(devoid of physics), whereas, in this work, the topic is solving a 
classical nonlinear boundary value problem in mechanics.
If the iterative stiffening algorithm
cannot untangle the mesh, then UBN reports failure to solve the 
problem. 
\item
Else if iterative stiffening succeeds, then we take Newton iterations to solve
\eref{eq:varform}.  The starting point for Newton is an untangled mesh
produced by step 1.  No continuation is used.  On the 
other hand, Newton's method is safeguarded using a line search described
in Section~\ref{sec:NewtonLineSearch}, which prevents the introduction
of new tangling. The line search is based on a technique common in the 
interior point literature (see e.g.,~\cite{wright_ipm}).
\end{enumerate}

\section{Newton continuation}
\label{sec:contin}

In the case that direct use of Newton's method to find $\phi$ fails to
converge, the standard alternative is Newton continuation, also known
as applying the load in steps.  This section
briefly describes Newton continuation before we return
to a description of UBN.

The basic form of Newton continuation is quite straightforward: a
sequence of parameters $0=\tau_0<\tau_1<\cdots<\tau_N=1$ is
chosen, and a sequence of displacement vectors $\u_0,\u_1,$\\
$\ldots,\u_N$
is computed, in which for each $k$, $\u_k$ is the solution to
the discretized
\eref{eq:varform} in the case that $(\phi_0(\XX)-\XX,\b,\t_0)$
are replaced by
$\tau_k\cdot(\phi_0(\XX)-\XX,\b,\t_0)$.  Solution $\u_0$ (corresponding
to absence of loads) is identically $\bz$.
(In the case that additional information is available about the
final solution, one might be able to formulate a better initial
guess for $\u$; however, $\u_0=\bz$ is the default value for
most continuation codes.)
Solution $\u_k$ is found via Newton's method, where
$\u_{k-1}$ is used as the initial guess.  
The final deformed configuration
is given by $\u_N$ since $\tau_N=1$.  Note also that it is possible
to accept a low-accuracy (not fully converged)
solution for $\u_k$ when $k<N$ since it is
presumably not necessary to achieve high accuracy for intermediate results
that are not part of the ultimate answer.

In some cases, a straight linear parametrization of the load path (as
in the previous paragraph) is not feasible.  In this case, one must
construct
a nonlinear parametrization $(\phi_{NLP}(\XX;\tau), \b_{NLP}(\XX;\tau),
\t_{0,NLP}(\XX;\tau))$ with the property that 
$\phi_{NLP}(\XX;0)=\XX$
while 
$\phi_{NLP}(\XX;1)=\phi_0(\XX)$
and similarly for the other load
terms.  Examples of nonlinear parametrizations are given later in the
paper.

The only remaining issue is how to select the sequence of $\tau_k$'s.
We use an adaptive rule defined as follows.  
Assume that there are no body forces and that the traction boundary 
conditions
are all zero (i.e., ``traction-free'' surfaces).  This means that the only
loading term is the Dirichlet boundary condition.
We form the deformed
mesh $M_{k-1}$ after applying the displacements given 
by $\u_{k-1}$ to non-Dirichlet nodes and
Dirichlet boundary conditions scaled by $\tau_{k-1}$,
(i.e., the deformed position is given by 
$\XX+\tau_{k-1}(\phi_0(\XX)-\XX)$
in the case of linear parametrization)
to Dirichlet nodes.
Next, we
compute a value of $\tau_k$ such that, if the boundary nodes 
in $M_{k-1}$  are further
deformed to positions given by $\XX+\tau_k(\phi_0(\XX)-\XX)$, 
then no tetrahedron
altitude will decrease by more than a factor of $\eta$, where
$\eta$ is a tuning parameter of the continuation algorithm.  Typically 
$\eta\le 1$.
(In particular, the step is sufficiently small that the mesh will not 
tangle after the
new boundary condition given by $\tau_k$ is applied to $M$.)
We also investigate some more aggressive continuation strategies with 
$\eta > 1$ in our 
experiments in Sections \ref{sec:2Dexperiments} and 
\ref{sec:3Dexperiments}.  
This adaptive strategy appears to work reasonably well, although we
did encounter some robustness problems discussed in 
Section~\ref{sec:3Dexperiments}.  We also compare these adaptive step 
selection strategies with a constant step-size strategy.

In this paper, we assume that the problem under consideration is to
determine a single final configuration.  Newton continuation 
finds this final configuration, and, as a by-product, 
also computes many intermediate configurations.  In some applications
this ``by-product'' is in fact the principal application of continuation.
For example, the
entire loading path is sometimes sought when the hyperelastic material
is, in and of itself, the object of study (e.g., a study of
soft-tissue deformation or damage due to an impact).  

On the other hand, for problems in which the hyperelastic material is
merely one component of a larger problem (e.g., a vibration isolator
in the model of a large structure), the entire load path is usually
not needed.  Furthermore, even in applications where the entire
loading path is required, our technique is applicable since UBN can be
used in combination with Newton continuation to obtain an improved
initial guess and larger steps than is possible using Newton
continuation alone.

The description in the earlier paragraphs assumed the special case 
of traction-free Neumann boundaries and absence of body forces.
It is more difficult to use this
adaptive technique when there are nonzero body forces or tractions since
it is not obvious
how to step these loads in a way that prevents tangling on each
step. Therefore, most of our test cases focus on the traction-free
case.
Since the focus of the paper is the UBN method, 
it represents a strengthening of
our contention that UBN is usually better than the
competing algorithm (continuation) since
we limit our testing only to the case
that seems well suited for 
continuation.
Nonetheless, we have also tried examples with nonzero tractions; we
report on this experiment at the end of Section~\ref{sec:3Dexperiments}.

\section{Iterative stiffening for mesh untangling}
\label{sec:femwarp}

In this section, we describe our procedure called iterative stiffening
for untangling a mesh.  We take the original mechanical problem
given by \eref{eq:varform}, and using
the same boundary conditions and loads, we solve the equations of isotropic
linear elasticity using piecewise linear (constant-strain) finite
elements 
\cite{BelytschkoLiu}.
Note that these equations have the same material
parameters (the Lam\'e constants $\lambda$ and $\mu$) as the Mooney-Rivlin
model.  Linear elasticity 
requires one linear system solve.  If the deformed mesh
(i.e., the mesh that arises from moving the nodes to their displaced
positions) is untangled, then iterative stiffening is finished.
If not, then our iterative stiffening procedure locates all elements
that are inverted in the deformed mesh and increases their stiffness
by 50\%.  The linear elasticity model is now solved again.
This procedure is repeated indefinitely until the mesh is
untangled or an excessive number of iterations has passed.

We have not found this precise version of
iterative stiffening  
appearing in the previous literature, but it is related to 
ideas already in the literature.
It is closely related to ``Jacobian techniques'' of
Stein et al.\ \cite{SteinTezduyarBenney}.  
It could be regarded as an extension of FEMWARP, a finite element based mesh warping approach 
developed by the authors within the linear weighted 
Laplacian smoothing (LWLS) framework~\cite{ShontzVava,shontz-thesis}.
One difference is that FEMWARP does not easily encompass the mechanical
concept of traction boundary conditions.
It is also related to a mesh warping method used for ALE
solvers and described in Chapter 7 of 
\cite{BelytschkoLiu}.

We remark that iterative stiffening, which we treat herein as the first
step of UBN, could be a standalone algorithm for ALE remeshing.
Indeed, this is the application for ``Jacobian techniques'' mentioned
above.  

In our preliminary version of the UBN method \cite{shontz-thesis}, 
the untangling was done using Opt-MS~\cite{FreitagPlassmann} rather
than iterative stiffening.  
Opt-MS is an untangling algorithm that
iteratively repositions interior nodes one at a time until
the mesh is untangled.  It solves a small linear-programming problem for 
each node to find the position for it that maximizes the minimum
area (volume) of an element in the local submesh constructed from its
neighboring triangles (tetrahedra).  The area (volume) of a triangle 
(tetrahedron) is computed via the determinant of the Jacobian of the 
element.  We found recently that 
iterative stiffening is more effective for use in UBN than Opt-MS.  One 
possible reason is that it is difficult
to implement traction boundary conditions in a natural way in Opt-MS.

Note that iterative stiffening can be made particularly efficient
by using matrix-updating.  In particular, it is well-known (see,
e.g. \cite{GVL}) how to update a Cholesky factorization of 
a symmetric positive definite matrix $A$ after $A$ has undergone
a low-rank update.  If the iterative stiffening procedure
stiffens only a few
elements per iteration (our test runs confirm that indeed there
are usually only a few updates per step), then this can be implemented as a 
a low-rank update, which is potentially much more efficient
than solving a new stiffness matrix from scratch.  We did not implement
matrix-updating 
because the work for iterative stiffening was usually dominated
by the solver part of the algorithm anyway.

\section{Newton Line Search}
\label{sec:NewtonLineSearch}
Newton's method is often employed for solving nonlinear systems of 
continuously differentiable equations~\cite{dennis_schnabel}.  Let 
$f \mbox{:} \R^{dm} \rightarrow \R^{dm}$, continuously differentiable, be given.  The 
task at hand is to find a $\u \in \R^{dm}$ such that $f(\u) = \bz$.  
Let $\u_0 \in \R^{dm}$ be given.  Then, at each iteration $k$, 
Newton's method solves \begin{equation}
J(\u_k)\s_k = -f(\u_k),
\label{eqn:newton}
\end{equation}
where $J$ denotes the Jacobian of $f$, for the Newton step, $\s_k$, and 
performs the following update
\begin{equation}
\u_{k+1} = \u_k+\s_k.  
\label{eqn:newton_update}
\end{equation}
If it becomes necessary to satisfy one or more additional inequality
constraints, it is 
possible to safeguard the Newton step with the introduction of a line
search.  Let $\alpha_k$ denote the line search parameter.  Then
$\alpha_k$ is chosen to be as large as possible such that 
$0 < \alpha_k \leq 1$ and $\u_{k+1} = \u_k +\alpha_k \s_k$ satisfies the 
constraint.  

It is often difficult to compute the value of $\alpha_k$ that 
minimizes $f(\x_k+\alpha_k \s_k)$ and satisfies the constraint
because $f$ is often a highly nonlinear function.  In addition, 
the optimal value of $\alpha_k$ often produces steplengths that are too 
short in practice.  Thus, it is common practice in interior point 
methods to derive heuristics for computing $\alpha_k$ that 
allow for both ease of computation and larger 
steplengths~\cite{wright_ipm}.  One such heuristic is to choose 
$\alpha_k$ so as to stay a fixed percentage away from the boundary.
We employ this heuristic in our line search below.

As was pointed out in Section~\ref{sec:tangling}, the mesh is tangled
unless $J = \det(F) > 0$.  Thus, we introduce a line search that enforces
that $J > 0$ on each iteration of Newton's method.  In particular, we
begin with $J > 0$ on the zeroth iteration and choose the line search
parameter $\alpha_k$ such that $J(\u_{k+1}) \geq 0.1 J(\u_k)$ on each
element so as to stay a fixed percentage away from the boundary for 
reasons discussed above.  

The following pseudocode algorithm 
shows how the line search parameter is determined.
Let $N$ denote the number of elements in the mesh.  
Given a displacement vector $\u$,
it is straightforward for each element $i=1,\ldots,N$
to compute the deformation gradient
$F$ and its determinant $J$
determined by this displacement on element $i$;
we denote the resulting determinant by $J(\u,i)$.

Let $\u_0$ be the value of the displacement (at non-Dirichlet nodes)
returned by the previous iteration of our Newton/line search
algorithm. Initially, the value of $\u_0$ is
the output of the iterative
stiffening algorithm.  It is assumed that the mesh 
determined by the Dirichlet boundary conditions and by
$\u_0$ on non-Dirichlet nodes is untangled.
Let $\s$ denote the Newton step determined from $\u_0$ via \eref{eqn:newton}.

\begin{tabbing}
++\=+++\=+++\=+++\kill
$\alpha = 1;$ \\
for $i = 1 \mbox{:} N$ \\
\>  while true \\
\> \>  if $J(\u_0+\alpha\s,i) \ge J(\u_0,i)/10$ \\
\> \> \> break \\
\> \> end \\
\> \> $\alpha = \alpha \cdot 0.9;$ \\
\> end \\
end
\end{tabbing}

\section{2D Experiments}
\label{sec:2Dexperiments}

We designed a series of numerical experiments
in order to test the robustness of UBN and to
compare it to the standard Newton continuation algorithm.  
As explained in Section \ref{sec:contin}, most of our 
test cases involve only traction-free, body-force-free loading conditions.
For all of the numerical experiments in this paper, we set the
parameters in \eref{eq:MR-model} as follows:  $\lambda = \frac{\nu 
E}{(1+\nu)(1-2 \nu)}$ and $\mu = \frac{E}{2(1+\nu)}$, with $E = 1$ and 
$\nu = 0.3$.  

The termination criteria for the Newton loop in UBN and for
the final step of Newton continuation was that
$\|F\|_2 \leq 10^{-10} \|F_0\|_2$, where $F_0$ is the initial
value (i.e., the value when all interior displacements are
set to 0) of the load vector.  For the 
Newton continuation steps prior to the final step,
the termination criteria was that
$\|F\|_2 \leq tol \|F_{k_i}\|_2,$ where
$F_{k_i}$ is the initial value of the load vector at the beginning
of major iteration $k$, and $tol = 10^{-3}$ or $10^{-5}$.
The looser tolerance was chosen because it was important to 
determine the value of the stopping criterion which makes Newton 
continuation as efficient as possible (for the purposes of comparison with 
UBN).  The tighter tolerance was chosen for the purposes of 
improving the robustness of Newton continuation on highly 
deformed meshes.  The algorithms were implemented in Matlab.



The linear solution
operation in Matlab is quite highly optimized and is expected to
compete well with a custom-written C or C++ linear solver.  On the
other hand, the matrix assembly process  involves several nested
Matlab loops and is therefore expected to be much slower than
a C or C++ version.  For this reason, wall-clock times derived from
the Matlab code are not useful predictors of computational demands
that would be observed with a C or C++ code.

Instead, we measure the running time in terms of assembly/linear solve
steps.  An assembly/linear solve (ALS) step 
consists of one stiffness matrix
and load vector assembly operation followed by one sparse linear
system solve.  The Newton continuation method involves a sequence
of Newton solve procedures, and each Newton solve is further subdivided
into several ALS steps.  The UBN method involves iterative stiffening
iterations 
followed by a safeguarded Newton method.  We count each iteration
of iterative stiffening as an ALS step.  The assembly portion of the
iterative stiffening ALS operation is not exactly the same as the assembly
portion of Newton, since the former involves linear elasticity assembly
whereas the latter involves nonlinear tangent stiffness assembly.  
We ran both assembly codes on an older Windows machine running Matlab 5.3,
which has a ``flops'' function built in that measures floating point 
operations.
(Newer versions of Matlab lack this function.)  From this experiment we
determined that the number of operations for the two kinds of assembly
are fairly close.  Furthermore, both assembly operations are much less
costly than the linear system solve.  Note that the iterations of
iterative stiffening would be considerably cheaper than an ALS step
had we implemented low-rank corrections described in 
Section~\ref{sec:femwarp}.

The solver portion of UBN involves additional operations connected
with the line search.  We determined (again by running test cases in
Matlab 5.3) that the line search requires a tiny number of operations
in comparison to the solution of the linear equations.  

Thus, it is sensible to compare the running time of
UBN to continuation by considering
the total number of ALS steps required for either.

\begin{figure}
\centering\includegraphics[height=5cm]{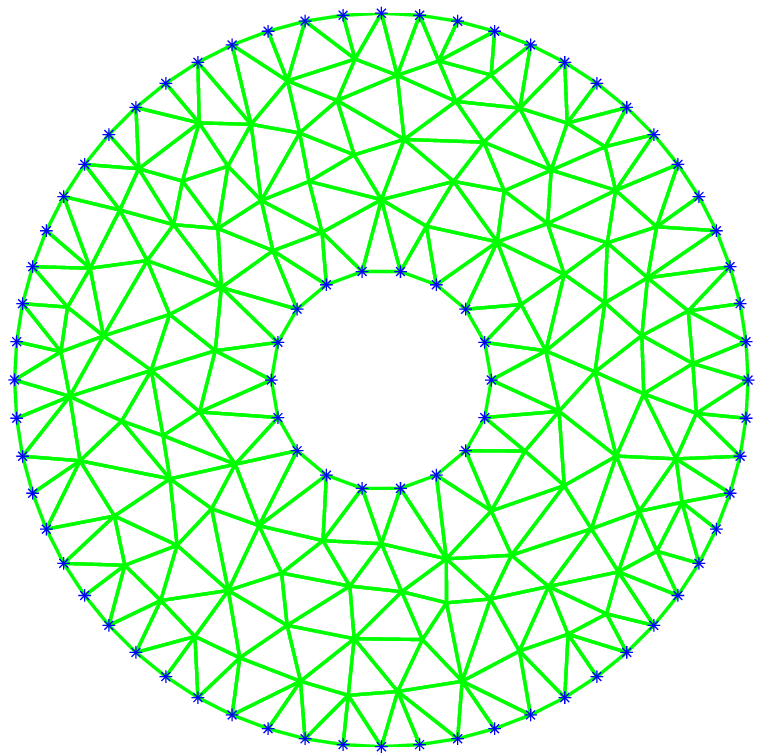}
\caption{The annulus mesh used for testing in this section.}
\label{fig:annmesh}
\end{figure}

In this section we describe our 2-dimensional test case, which
is an annular domain.  The mesh was generated with Shewchuk's
Triangle \cite{triangle} and is illustrated in Fig.~\ref{fig:annmesh}.
It contains 181 nodes and 284 triangles.

\begin{figure}
$$
\begin{array}{cc}
\centering\includegraphics[height=4cm]{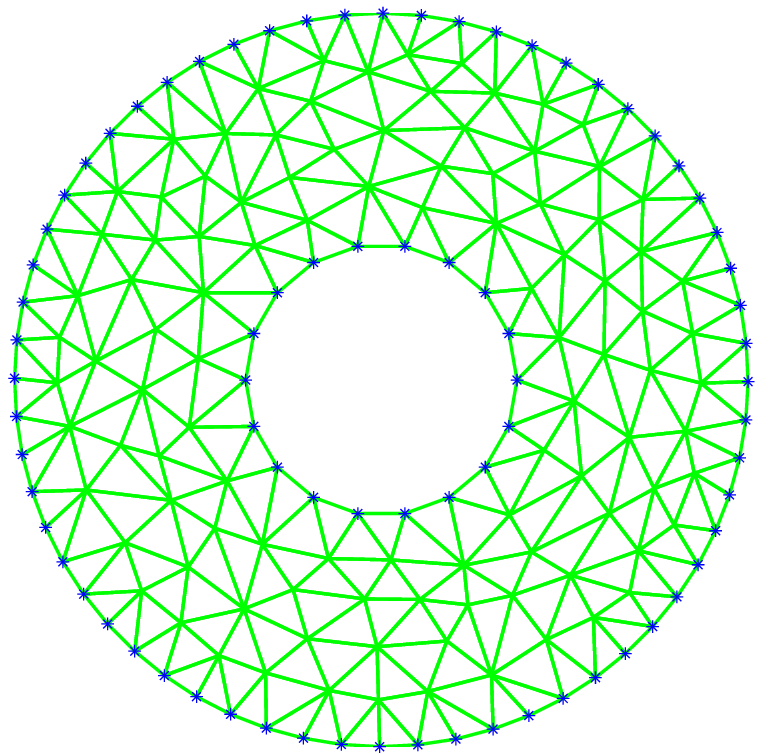} &
          \includegraphics[height=4cm]{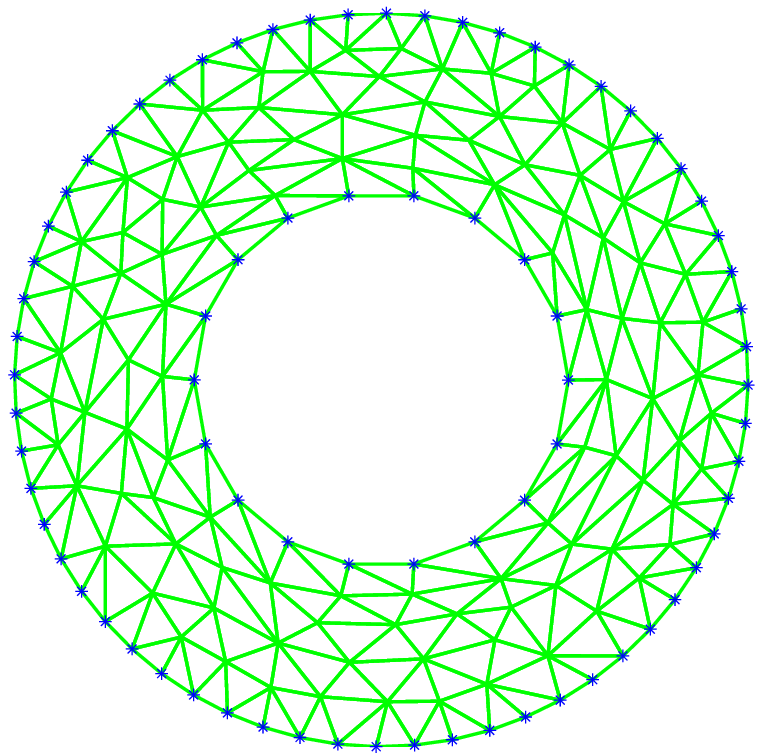} \\
          (a) & (b) \\
          \includegraphics[height=4cm]{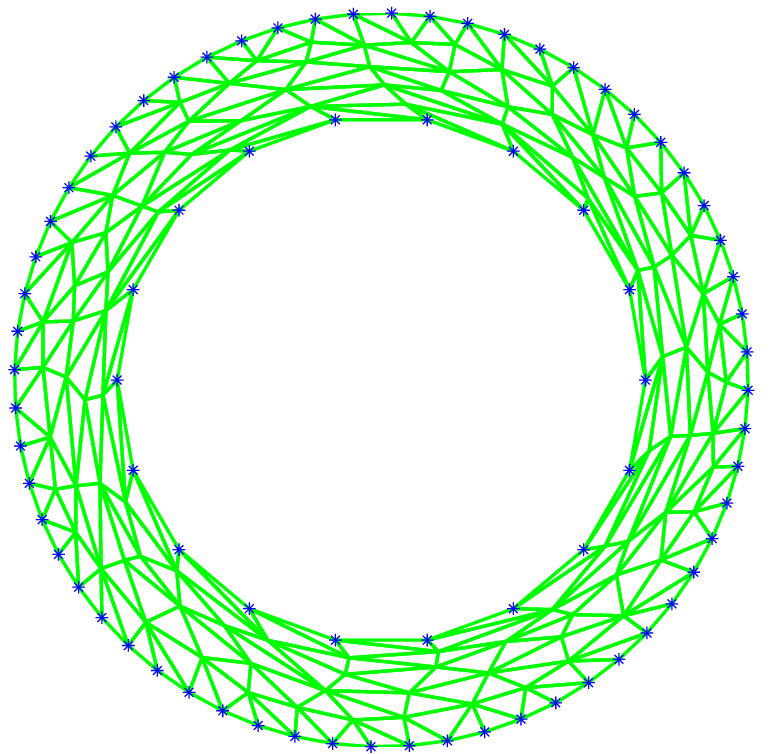} &
          \includegraphics[height=4cm]{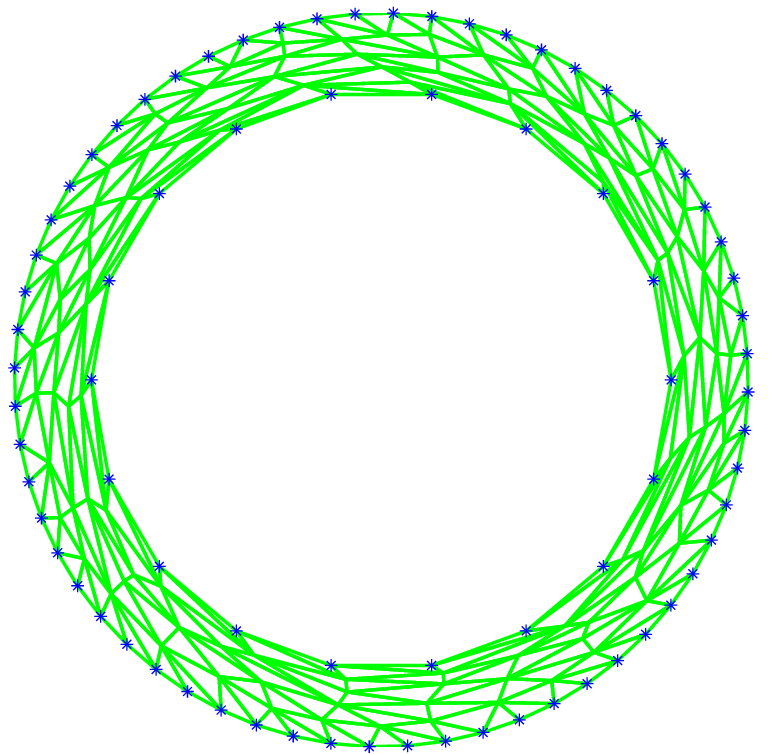} \\
          (c) & (d)
\end{array}
$$
\caption{Deformed annulus meshes resulting from rotating the exterior 
boundary circle of the mesh shown in Fig.~\ref{fig:annmesh} by $f$ radians 
and moving the inner boundary circle by a factor $f$ closer to the outer 
boundary.  The deformed meshes are for (a) $f=0.1$, (b) $f=0.3$, 
(c) $f=0.6$, and (d) $f=0.7$.}
\label{fig:anndef}
\end{figure}

The boundary conditions used in this test case involve a rotation
of the exterior boundary circle by $f$ radians combined with 
moving the inner boundary by a factor $f$ closer to the outer boundary
(where $f=0$ means no motion and $f=1$ means that the inner boundary 
would coincide with the outer boundary).  Values of $f$ tried were 
$0.1$, $0.3$, $0.6$, and $0.7$.  The resulting deformed meshes are 
illustrated in Fig.~\ref{fig:anndef}.  

The number of ALS steps to compute
these deformed meshes is given in Table~\ref{tab:annruns}.
The columns of this table are as follows.  The first column
$f$ is the amount of boundary deformation as described in the previous
paragraph.  The second column \# inv is the number of inverted elements 
in the deformed mesh prior to application of UBN.
The third column UBN--IS is the number of iterations of iterative
stiffening required by UBN.   The fourth column UBN--NM is the number
of iterations of Newton's method required by UBN.  
The fifth column UBN-ALS is the number of ALS steps required by UBN (and 
hence is the sum of the second and third columns).  

The remaining columns of the table report on results from the
continuation algorithm.  The sixth column is the number of major 
iterations (i.e., updates to the continuation parameter $\tau$) required 
by the continuation algorithm when constant-size steps are employed.
The seventh column is the number of ALS steps required by continuation.
The eighth and ninth columns are the same quantities required by the 
continuation algorithm using the adaptive rule discussed in 
Section~\ref{sec:contin} with parameter $\eta=1/3$.
The tenth and eleventh columns are the same quantities when
$\eta=1.2$.  Note that $\eta=1.2$ is a quite aggressive choice of
stepsize for continuation, since any value of $\eta>1$ means that 
updating the boundary could cause an inversion.  For this 2D test case,
an aggressive choice of $\eta$ did not seem to hinder convergence,
but the results for a large value of $\eta$ in 3D described in the
next section are less favorable.

For $f = 0.8$, neither UBN nor Newton continuation was able to find 
a solution. UBN's iterative stiffening did not untangle the mesh after the 
maximum number of iterations (400) had been reached, and continuation 
stalled at $\tau = 0.93$ when adaptive steps were used and terminated 
with an inverted element when constant steps were used.

It should be noted that a highly deformed mesh like the solution when 
$f=0.7$ is probably not physically valid because the finite element 
discretization is no longer an accurate approximation to the underlying 
PDE.  Nonetheless, we include extreme cases like this because it is 
interesting to compare the two algorithms in limiting cases.  The test 
results show that UBN is much faster than continuation for both modest and 
extreme deformations.

Note that for continuation, the outer boundary motion (i.e., the Dirichlet 
boundary condition) is parametrized in polar coordinates by $\theta$, the 
rotation angle.  Linear parametrization with respect to the 
rectangular 
coordinates, $(x,y)$, would work poorly in this case because a linear 
deformation from the initial position of the outer boundary to the final 
position would cause the outer boundary to shrink in radius and then 
expand.  

\begin{table*}
\caption{Results of comparison of UBN to the most efficient version of continuation for the 2D annular domain.  
See the corresponding text for explanation of the column headers.}
\label{tab:annruns}
\begin{center}
\begin{small}
\begin{tabular}{rrrrrrrrrrr}
\hline\noalign{\smallskip}
\multicolumn{1}{c}{$f$}& \multicolumn{4}{c}{UBN} &
\multicolumn{6}{c}{Contin., tol = $10^{-3}$} \\
& & & & & \multicolumn{2}{c}{const. steps}
& \multicolumn{2}{c}{$\eta=1/3$}
& \multicolumn{2}{c}{$\eta=1.2$} \\
& \multicolumn{1}{c}{\# inv} 
& \multicolumn{1}{c}{IS} & 
\multicolumn{1}{c}{NM} & 
\multicolumn{1}{c}{ALS} & 
\multicolumn{1}{c}{MajIt} &
\multicolumn{1}{c}{ALS} &
\multicolumn{1}{c}{MajIt} &
\multicolumn{1}{c}{ALS} &
\multicolumn{1}{c}{MajIt} &
\multicolumn{1}{c}{ALS}
\\
\hline\noalign{\smallskip}
0.1 & 0 & 1 & 3 & {\bf{4}} & 27 & 55 & 28 & 57 & 8 & {\bf{18}}  \\
0.3 & 36 & 1 & 5 & {\bf{6}} & 80 & 162 & 87 & 175 & 24 & {\bf{50}} \\
0.6 & 59 & 5 & 29 & {\bf{34}} & 160 & 322 & 265 & 483 & 73 & {\bf{148}} \\
0.7 & 64 & 9 & 23 & {\bf{32}} & 187 & 356 & 414 & 683 & 113 & {\bf{228}} \\	
\hline\noalign{\smallskip}
\end{tabular}
\end{small}
\end{center}
\end{table*}

Comparing the UBN--ALS and Contin--ALS columns of this
table indicates that UBN is approximately 9-27 times
more efficient than continuation when constant-size steps
are used, and is 15 to 30 times more efficient when $\eta=1/3$.  (Note 
that the most efficient results for the UBN and Newton continuation 
methods are shown in bold face type.)
Continuation is $3$ to $4$ times faster when used with the larger value of 
$\eta$ but is still significantly slower than UBN.
Other annulus
deformation tests not reported here confirm that UBN is always
far more efficient than continuation.

We also wished to check whether the iterative stiffening step in UBN was 
essential.  As evidenced in Column 2 of Table~\ref{tab:annruns}, for 
smaller deformations, the deformed mesh does not always result in inverted 
elements.  However, for larger deformations, the deformation does result 
in inverted elements. For deformed meshes with inverted elements, the 
iterative stiffening step is essential to untangling the mesh before using 
it as a starting point to the line search.  For all deformed meshes, it is 
useful for determining a good starting point.

Similarly, we checked to determine whether the line sea-\\rch procedure 
built 
into UBN was ever active in order to determine whether it is an essential
part of UBN.  We found that it was active on about 30\%-50\% of the
iterations for the larger values of deformation.

\begin{table*}
\caption{Results of comparison of different versions of continuation for the 2D annular domain.  
See the corresponding text for explanation of the column headers.}
\label{tab:annruns2}
\begin{center}
\begin{small}
\begin{tabular}{rrrrrrrrrrrrrrrrr}
\hline\noalign{\smallskip}
\multicolumn{1}{c}{$f$}& 
\multicolumn{10}{c}{Contin., tol = $10^{-3}$} &
\multicolumn{6}{c}{Contin., tol = $10^{-5}$} \\
& \multicolumn{2}{c}{const. steps}
& \multicolumn{2}{c}{$\eta=1/3$}
& \multicolumn{2}{c}{$\eta=1.2$} 
& \multicolumn{2}{c}{$\eta=1/3$ + LS} 
& \multicolumn{2}{c}{$\eta=1.2$ + LS}
& \multicolumn{2}{c}{const. steps} 
& \multicolumn{2}{c}{$\eta=1/3$}
& \multicolumn{2}{c}{$\eta=1.2$} \\
& \multicolumn{1}{c}{MajIt} &
\multicolumn{1}{c}{ALS} &
\multicolumn{1}{c}{MajIt} &
\multicolumn{1}{c}{ALS} & 
\multicolumn{1}{c}{MajIt} &
\multicolumn{1}{c}{ALS} &
\multicolumn{1}{c}{MajIt} &
\multicolumn{1}{c}{ALS} &
\multicolumn{1}{c}{MajIt} &
\multicolumn{1}{c}{ALS} &
\multicolumn{1}{c}{MajIt} &
\multicolumn{1}{c}{ALS} &
\multicolumn{1}{c}{MajIt} & 
\multicolumn{1}{c}{ALS} &
\multicolumn{1}{c}{MajIt} &
\multicolumn{1}{c}{ALS}
\\
\hline\noalign{\smallskip}
0.1 & 27 & 55 & 28 & 57 & 8 & {\bf{18}} & 28 & 57 & 8 & 18 & 27 & 55 & 28 
& 57 & 8 & 25 \\
0.3 & 80 & 162 & 87 & 175 & 24 & {\bf{50}} & 87 & 175 & 24 & 50 & 80 & 190 
& 87 & 
179 & 24 & 73 \\
0.6 & 160 & 322 & 265 & 483 & 73 & {\bf{148}} & 265 & 483 & 73 & 148 & 160 
& 430 & 264 & 545 & 73 
& 220  \\
0.7 & 187 & 356 & 414 & 683 & 113 & {\bf{228}} & 414 & 683 & 113 & 228 
& 187 & 511 & 407 & 832 & 
113 & 340 \\	
\hline\noalign{\smallskip}
\end{tabular}
\end{small}
\end{center}
\end{table*}

Motivated by the robustness issues experienced by the Newton continuation 
method in 3D, we also wished to investigate the performance of the method
when safeguarded steps are taken (by employing the same line search 
described in Section~\ref{sec:NewtonLineSearch}), a tighter convergence 
tolerance is used, and constant-size steps are employed.  The results of 
our experiments are shown in Table~\ref{tab:annruns2}.  For ease of 
comparison, the Newton continuation results shown in 
Table~\ref{tab:annruns} are listed here in columns 2-7.  The eighth and 
ninth columns record the results of adding a line search to the Newton 
continuation method for the case when $\eta=1/3.$  The tenth and eleventh 
columns indicate the corresponding information when a line search is used 
and $\eta=1.2.$  The twelfth and thirteenth columns show the results when 
constant steps are employed with a convergence tolerance of $10^{-5}.$  
The remaining columns specify the results for when the 
convergence tolerance is $10^{-5}$ and $\eta=1/3$ or $\eta=1.2.$

The results show that the addition of a line search to safeguard the steps 
of the Newton continuation method is unnecessary whenever the method is 
able to compute the deformed mesh without experiencing robustness issues.  
However, the 
continuation method was able to compute the deformed mesh for $f=0.8$ when 
a convergence tolerance of $10^{-3}$ and a line search were employed.
The use of a tighter convergence tolerance slows down the 
convergence of the continuation method on these test problems.  Taking 
constant-size steps generally produces convergence more quickly than the 
use of the conservative adaptive stepping rule but more slowly than the 
aggressive adaptive stepping rule.  However, it should be noted that the 
use of a convergence tolerance of $10^{-5}$ when taking 
either constant-size steps or steps computed using the aggressive adaptive 
step 
rule were also able to compute a deformed mesh for $f=0.8$.  In addition, 
the faster Newton continuation methods discussed above are up to 4 times 
faster than the other continuation methods presented here.  Thus, 
the use of a 
line search and a narrower convergence tolerance should be performed only 
when the deformation is extremely large, and UBN is not able to compute 
the deformed mesh.

The use of a line search with a convergence tolerance of $10^{-5}$ was not 
explored since the use of constant-size steps with the tighter 
convergence tolerance was the most effective (i.e., the method was even 
able to compute the deformed mesh for $f=0.9$).

\section{3D Experiments}
\label{sec:3Dexperiments}

Our experiments in 3D consisted of two tetrahedral meshes called ``Hook''
and ``Foam5,'' which were provided to us by P.~Knupp \cite{pat_mesh}.
``Hook'' is a geometry composed of three main sections:  its two end 
segments are composed of half annuli (in 3D), and its middle section 
is an irregularly-shaped solid which creates a sharp corner where it 
joins the bottom section.  ``Foam5'' is a prism whose cross-section is a 
half-disk with three cavities cut on the top surface; two of the cavities 
are cylinders and the third is two parallelpipeds arranged like stairs.
The sizes of the meshes are as follows: Hook contains 1190 nodes and 4675 
tetrahedra, and Foam5 contains 1337 nodes and 4847 tetrahedra.  Hook is 
contained in a bounding box of size $54\times 40\times 95$, while Foam5 is 
contained in a bounding box of size $11.3\times 5.5\times 6.6$.

\begin{figure}
$$
\begin{array}{cc}
\centering\includegraphics[height=5cm]{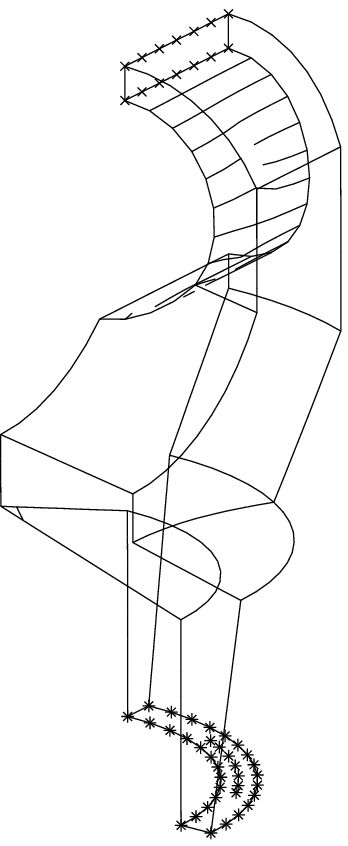} &
          \includegraphics[height=5cm]{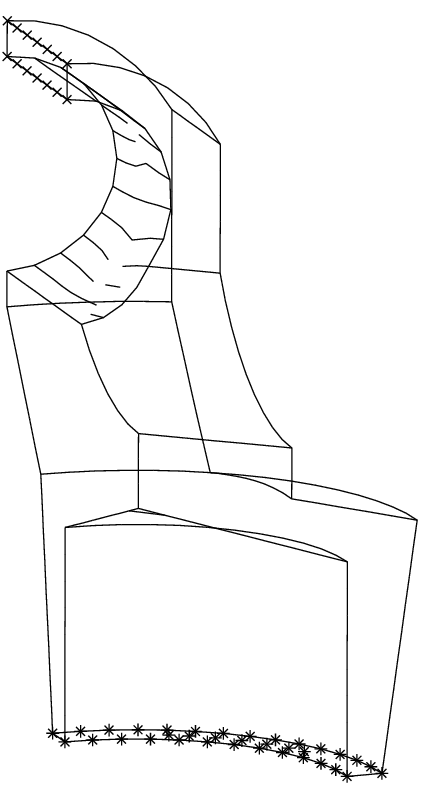} \\
          (a) & (b) \\
          \includegraphics[height=5cm]{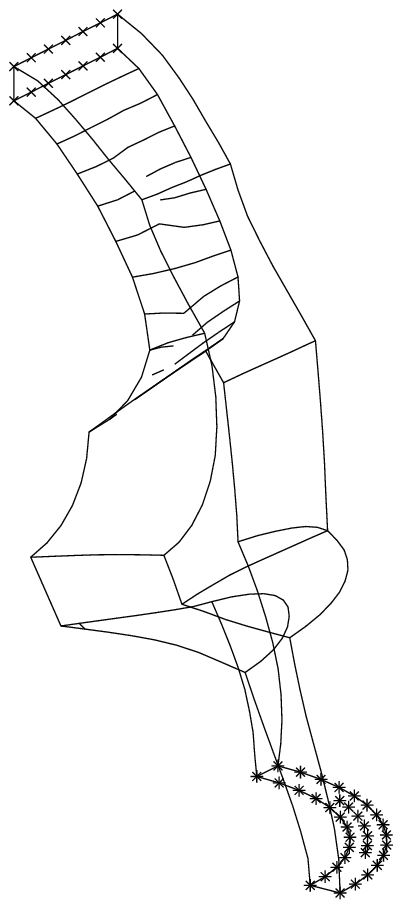} &
          \includegraphics[height=5cm]{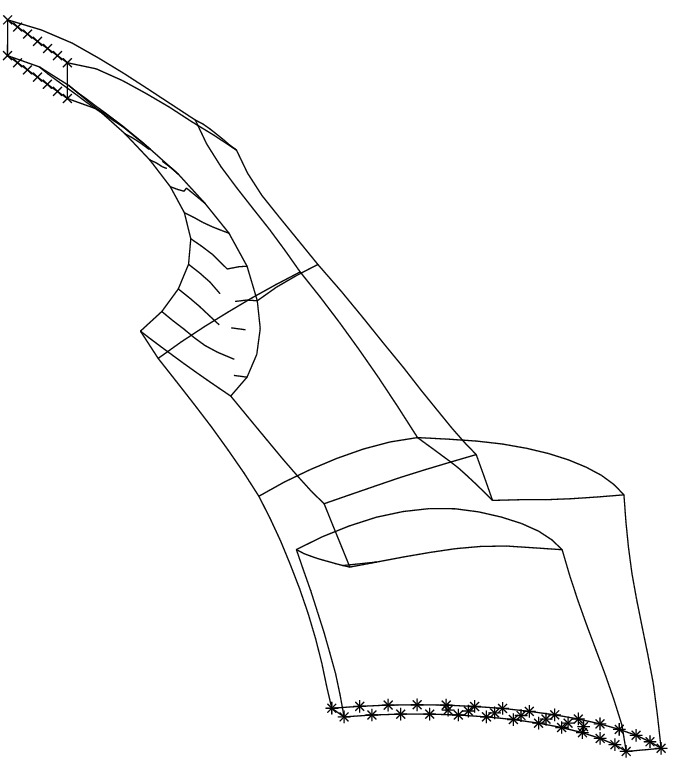} \\
          (c) & (d)
\end{array}
$$
\caption{The top row line diagrams $(a)$ and $(b)$
show the undeformed Hook body from two different
viewpoints.  Dirichlet boundary conditions were
applied to two of the boundary surfaces to yield
deformed meshes.  In particular, the asterisks mark the 
zero-displacement boundary, while the $\times$'s mark 
fixed displacement.
The bottom row diagrams $(c)$ and $(d)$ show Hook after the maximum
deformation of 40 is applied.}
\label{fig:hookfig}
\end{figure}

\begin{figure*}
$$
\begin{array}{cc}
\centering\includegraphics[height=5cm]{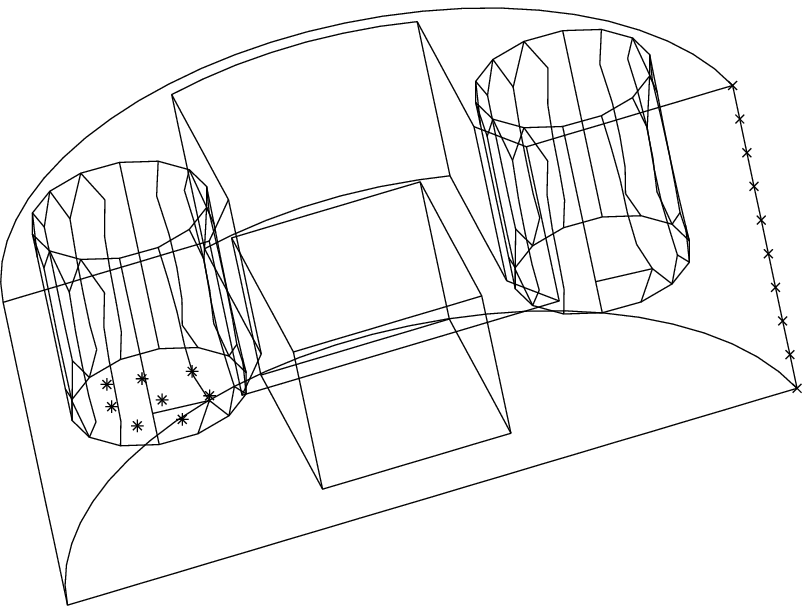} &
          \includegraphics[height=5cm]{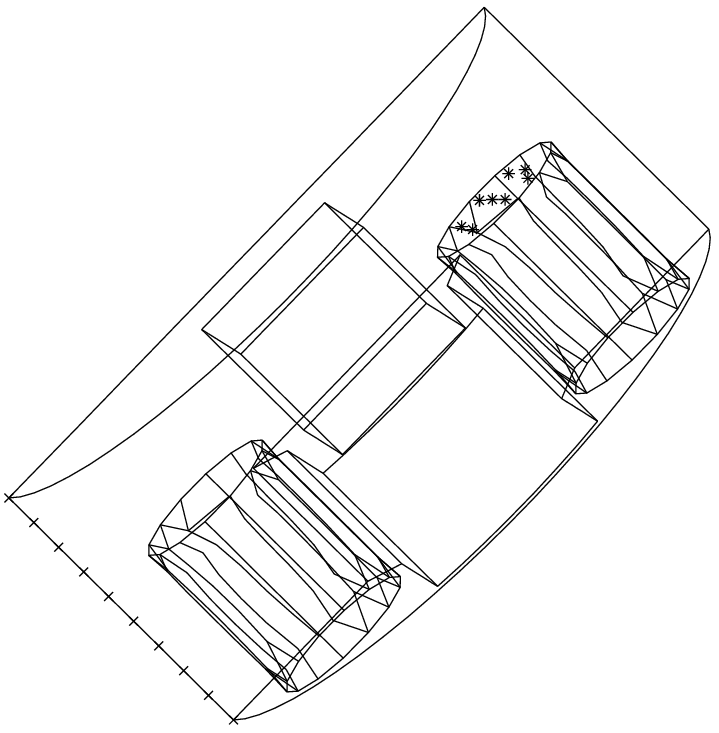} \\
          (a) & (b) \\
          \includegraphics[height=5cm]{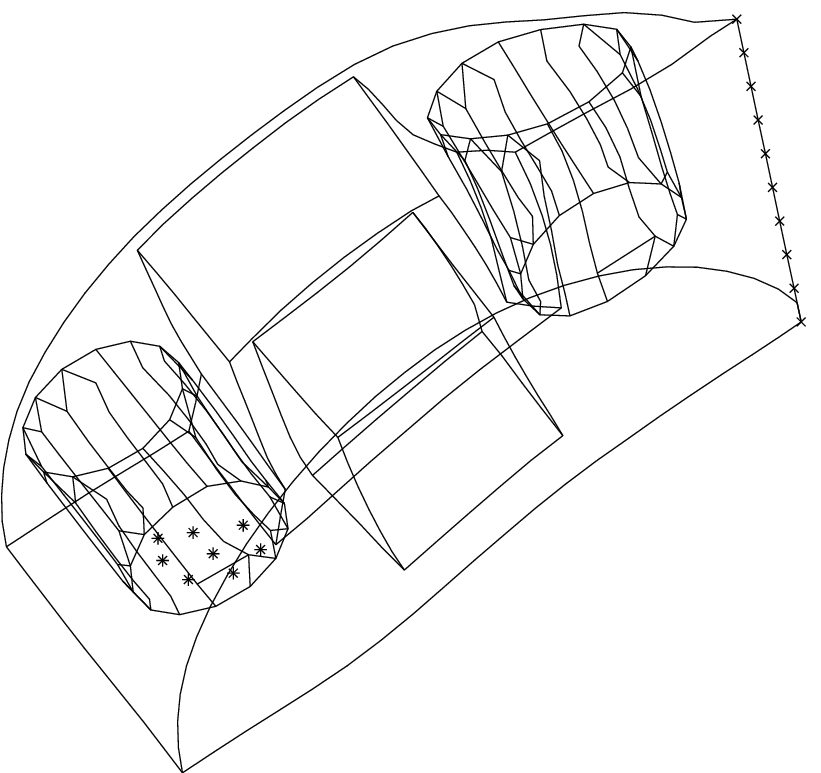} &
          \includegraphics[height=5cm]{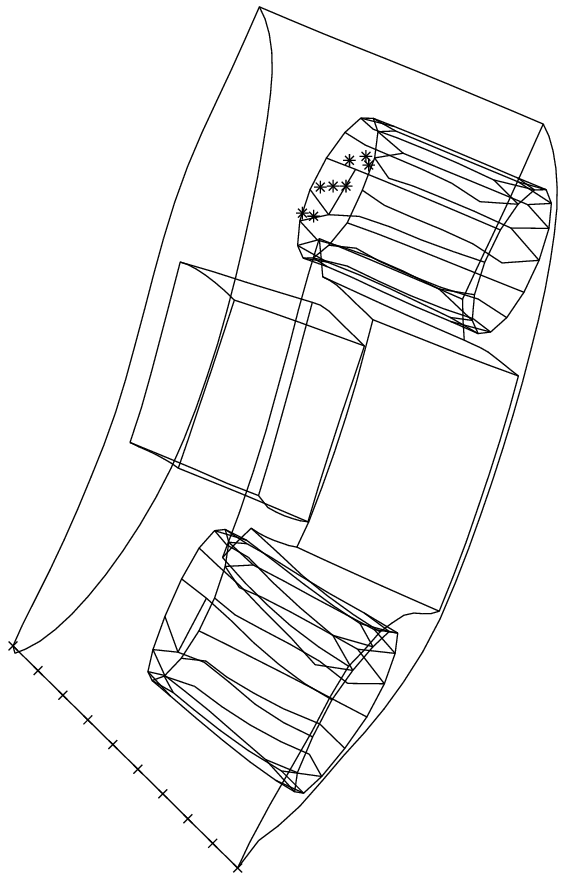} \\
          (c) & (d) 
\end{array}
$$
\caption{The top row line diagrams $(a)$ and $(b)$
show the undeformed Foam5 body from two different
viewpoints.  Dirichlet boundary conditions were
applied to two of the boundary surfaces to yield
deformed meshes.  In particular, the asterisks mark the 
zero-displacement boundary,
while the $\times$'s mark fixed displacement.
The bottom row diagrams $(c)$ and $(d)$ show Foam5 after the maximum
deformation of 5 is applied.}
\label{fig:foamfig}
\end{figure*}

In both cases, we applied Dirichlet boundary conditions to two of the
boundary surfaces, leaving the rest traction-free.  
In both cases, the Dirichlet
conditions are identically zero on one boundary
surface and displace the other surface in a uniform
direction.  Three magnitudes for the displacement were tested.
For Hook, the displacement sizes were
10, 20, and 40, whereas for Foam5 they were 0.5, 2, and 5.  Thus, we see 
that
the applied displacements are on the same order as the size of the
object, and therefore large deformations will result.
Figure \ref{fig:hookfig} shows the deformed and undeformed configurations
of Hook for the maximum deformation of 5, while Fig.~\ref{fig:foamfig}
shows the corresponding illustration of Foam5.

The results of our tests of UBN versus Newton continuation
on Hook and Foam5 with a $10^{-3}$ convergence tolerance are given in
Table~\ref{tab:3dresults}.  (The relevant column headers are the same as 
those given for Table~\ref{tab:annruns}.)  
As in the previous section, the unit
for measuring running time is ALS steps.  For these tests, straight
linear parametrization was used for continuation.
As before, our results on Hook indicate that UBN is 15 to 50 times
faster than continuation when $\eta=1/3$, and 5 to 15 times faster for
$\eta=1.2$.  In addition, UBN is 15 to 50 times faster than continuation
when constant-size steps are taken.   

The continuation algorithm terminates when 
the increment in $\tau$ becomes smaller than
the prespecified minimum (0.0005); this happened in two of
the tests with Foam5 when $\eta=1/3$ (indicated by `---'
in the table).
Apparently this is due to an extremely
flat tetrahedron which, although it is does not become inverted, causes
the heuristic used for adaptively incrementing $\tau$ to take very
conservative steps.  The continuation algorithm also terminates when 
inverted elements remain after one major iteration.  This happened in one 
of the tests with Hook when constant-size steps were taken as indicated by 
'***' in the table.  

For the Foam5 tests, the continuation algorithm terminated after one major 
iteration due to the presence of inverted elements for three of the meshes 
when $\eta=1.2$.  This shows that, as expected, the aggressive choice of 
$\eta$ may be more prone to inverting elements.  Similar performance of 
the continuation algorithm occurred for two of the meshes when 
constant-size steps were used.

\begin{table*}
\caption{Results of comparison of UBN to continuation for 3D domains.
See the corresponding text for explanation of the column headers.}
\label{tab:3dresults}
\begin{center}
\begin{small}
\begin{tabular}{lrrrrrrrrrrr}
\hline\noalign{\smallskip}
\multicolumn{1}{l}{mesh} & 
\multicolumn{1}{c}{displ.} & 
\multicolumn{4}{c}{UBN} &
\multicolumn{6}{c}{Contin., tol = $10^{-3}$} \\
& & & & & & \multicolumn{2}{c}{const. steps} &
\multicolumn{2}{c}{$\eta=1/3$} &
\multicolumn{2}{c}{$\eta=1.2$} \\
& & \multicolumn{1}{c}{\# inv}
& \multicolumn{1}{c}{IS} & 
\multicolumn{1}{c}{NM} & 
\multicolumn{1}{c}{ALS} & 
\multicolumn{1}{c}{MajIt} &
\multicolumn{1}{c}{ALS}  &
\multicolumn{1}{c}{MajIt} &
\multicolumn{1}{c}{ALS} &
\multicolumn{1}{c}{MajIt} &
\multicolumn{1}{c}{ALS} \\
\hline\noalign{\smallskip}
Hook & 10.0 & 0 & 1 & 5 & {\bf{6}} & 53 & 107 & 53 & 107 & 15 & 
{\bf{32}} \\
Hook & 20.0 & 0 & 1 & 6 & {\bf{7}} & 165 & 212 & 105 & 212 & 30 & 
{\bf{61}} \\
Hook & 40.0 & 0 & 1 & 8 & {\bf{9}} & {***} & {***} & 210 & 421 & 59 & 
{\bf{119}} \\
Foam5 & 0.5 & 18 & 1 & 4 & {\bf{5}} & {28} & {36} & 28 &  36  & *** & 
*** \\
Foam5 & 2.0 & 72 & 1 & 5 & {\bf{6}} & {***} & {***} & --- &  --- & *** & 
*** \\
Foam5 & 5.0 & 76 & 1 & 7 & {\bf{8}} & {***} & {***} & --- &  --- & *** & 
*** \\
\hline\noalign{\smallskip}
\end{tabular}
\end{small}
\end{center}
\end{table*}

\begin{table*}
\caption{Results of comparison of different versions of continuation for 3D domains.
See the corresponding text for explanation of the column headers.}
\label{tab:3dresults2}
\begin{center}
\begin{small}
\begin{tabular}{lrrrrrrrrrrrrrrr}
\hline\noalign{\smallskip}
\multicolumn{1}{l}{mesh} & 
\multicolumn{1}{c}{displ.} & 
\multicolumn{8}{c}{Contin., tol = $10^{-3}$} &
\multicolumn{6}{c}{Contin., tol = $10^{-5}$} \\
& & 
\multicolumn{2}{c}{$\eta=1/3$} &
\multicolumn{2}{c}{$\eta=1.2$} &
\multicolumn{2}{c}{$\eta=1/3$ + LS} &
\multicolumn{2}{c}{$\eta=1.2$ + LS} &
\multicolumn{2}{c}{const. steps} &
\multicolumn{2}{c}{$\eta=1/3$} &
\multicolumn{2}{c}{$\eta=1.2$} \\
& & \multicolumn{1}{c}{MajIt} &
\multicolumn{1}{c}{ALS}  &
\multicolumn{1}{c}{MajIt} &
\multicolumn{1}{c}{ALS} &
\multicolumn{1}{c}{MajIt} &
\multicolumn{1}{c}{ALS} &
\multicolumn{1}{c}{MajIt} &
\multicolumn{1}{c}{ALS} &
\multicolumn{1}{c}{MajIt} &
\multicolumn{1}{c}{ALS} &
\multicolumn{1}{c}{MajIt} &
\multicolumn{1}{c}{ALS} &
\multicolumn{1}{c}{Majit} &
\multicolumn{1}{r}{ALS} \\
\hline\noalign{\smallskip}
Hook & 10.0 & 53 & 107 & 15 & {\bf{32}} & 53 & 107 & 15 & 
32 & 53 & 159 & 53 & 
159 & 15 & 46 \\
Hook & 20.0 & 105 & 212 & 30 & {\bf{61}} & 105 & 212 & 30 
& 61 & 105 & 316 & 105 & 
316 & 30 & 91 \\
Hook & 40.0& 210 & 421 & 59 & {\bf{119}} & 210 & 421 & 
59 & 119 & xxx & xxx & 59 & 
177 & xxx & xxx\\
Foam5 & 0.5 & 28 &  36  & *** & *** & 28 & 36 & *** & 
*** 
& 28 & 57 & 8 & {\bf{25}} & 28 
& 57 \\
Foam5 & 2.0 & --- &  --- & *** & *** & 
--- & --- & *** & *** & 
110 & 221 &
109 & 219 & 31 & {\bf{93}} \\
Foam5 & 5.0 & --- &  --- & *** & *** & 
--- & --- & *** & *** & 274 & 549 & 269 & 539 
& 75 & {\bf{226}} \\
\hline\noalign{\smallskip}
\end{tabular}
\end{small}
\end{center}
\end{table*}

In order to attempt to improve the robustness of the Newton continuation 
algorithm for the 3D Hook and Foam5 mes-\\hes, we also performed 
experiments 
which considered the use of a line search to safeguard the steps, the use 
of 
constant steps, and the use of a tighter convergence tolerance.  The 
results of our experiments are shown in Table~\ref{tab:3dresults2}.  Note 
that the column headers are identical to those described above for 
Table~\ref{tab:annruns2} (except that the results for taking 
constant-size steps with a convergence tolerance of $10^{-3}$ has been 
omitted from the table).  The most efficient continuation method results 
are shown in bold face type.

The results demonstrate that the addition of a line search to the 
continuation method did not serve to resolve the robustness issues for 
the method when a convergence tolerance of $10^{-3}$ was employed.
However, the use of a tighter convergence tolerance, i.e., $10^{-5}$, 
(either in combination wi-\\th adaptive steps or constant-size steps) 
served 
to resolve the robustness problems seen when the looser convergence 
tolerance was employed.  The disadvantage is that the use of a tighter 
convergence tolerance makes the Newton continuation method much more 
expensive.  In 
particular, UBN is up to 70 times faster than the slowest Newton 
continuation method reported here.  It should be noted that for one of the 
Hook test cases, the continuation method terminated after the maximum 
number of iterations (600) had been performed; this was recorded as an 
'xxx' in the table.  The most 
efficient Newton continuation method is a function of the test problem.  
For somewhat smaller deformations (as was the case for the Hook mesh), the 
use of the aggressive adaptive step strategy and the looser convergence 
tolerance was the most efficient.  However, for larger deformations, 
the use of the aggressive adaptive step strategy with the tighter 
convergence 
tolerance was the most efficient.  Finally, the UBN method was much more 
efficient and more robust than the Newton continuation method, in general.

All the preceding tests involved traction free boundaries and
prescribed displacements, not all zero, for Dirichlet bou-\\ndary nodes.
We conclude this section by reporting on experiments
with the following boundary conditions.  Nonzero
tractions were specified on one facet of the Hook mesh, 
while zero displacements
were forced on a different facet.  (Tractions were implemented as
normally directed point lo-\\ads on each node of the facet.
For larger loads, this is not a completely realistic approach
since realistic forces would change directions under very large
deformation, but so-called ``follower'' loads are beyond the scope
of this work.)
The remaining boundary nodes
were traction-free.  Different levels of the traction load
were used for different experiment.

These experiments required a modification of our stepping rule for
continuation since the rule outlined in Section~\ref{sec:contin} is
intended for nonzero boundary displacements and zero tractions.
The modified continuation routine determines a fixed stepsize 
as follows.  First, the same underlying problem is solved using
linear elasticity.  (It may happen that some elements are inverted
in this solution; in this setting, we do not care about element
inversion.)  From this linear solution, we measure the maximum displacement
among nodes.  Then the stepsize for continuation is taken to be the
quotient of the minimum altitude in the original mesh divided by
the maximum displacement in the preliminary solve.  The rationale for
this rule
is so that the amount of deformation that occurs per step of continuation
should not exceed the sizes of the elements in an effort to prevent
inversions.

This stepsize rule appeared in our experiments to be appropriate in the
following sense.  Most outer iterations of continuation (i.e., 
stepping from $\tau_k$ to $\tau_{k+1}$) appeared to require 2 to 4
inner Newton iterations.  If the usual number required were 1, this
would indicate a stepsize which is too small (conservative).  
On the other hand,
if the usual required were much greater than 1, this would indicate
that the stepsize is too large for straightforward continuation. 

We found that UBN was 2 to 5 times faster than continuation
for these test cases.  Both algorithms returned a converged solution.
In the case of the largest load, the two solutions differed.  Both
were physically valid; one corresponded to the base of the hook
bending toward the hook end in the direction on the inside of the
hook, whereas the other corresponded to bending toward the outside
of the hook.  See remarks on the possibility of multiple solutions
in Section~\ref{sec:conclusions}.

Since all displacement boundary conditions in this example are zero,
the possibility of some additional experiments to elucidate features
of UBN were carried out.  The first experiment on this problem
ran the safeguarded Newton method of UBN but omitted the preliminary
use of FEMW-\\ARP to find a good starting point for the safeguarded 
Newton method.
Instead, the safeguarded Newton method was initialized with the original
mesh, which is possible because the prescribed displacement boundary
conditions are all zero.  We found that the method did not always
converge.  This shows that even the safeguarded Newton method should
be initialized close to the solution else divergence may result.
The second additional experiment looked at using unsafeguarded
Newton's method from the initial mesh to find the final configuration.
Again, this is possible because of the zero displacement condition.
Our experiment indicated that this method did not always converge
either.
Thus, these experiments provide evidence of the
necessity of the iterative stiffening and safeguarded line search.

\section{Conclusions}
\label{sec:conclusions}

In summary, we developed a robust solution method for solving
nonlinear elasticity equations for hyperelastic solids with large
boundary deformations.  The basic idea is to first untangle the mesh
using purely geometric methods and second solve the mechanical model;
thus, the algorithm was named UBN (for ``untangling before Newton'').
The first step of our algorithm is to attempt to untangle the mesh
with iterative stiffening. 
Assuming the mesh is untangled, UBN
takes safeguarded Newton steps to solve \eref{eq:varform}.

We tested the robustness of UBN and compared it to the standard Newton
continuation algorithm.  We demonstrated that UBN is significantly
more robust that the Newton continuation algorithm, i.e., it is
able to tolerate much larger deformations, in general.  For a couple of 
cases with extremely large deformations, the Newton continuation algorithm  
with the use of a tight convergence tolerance and very small 
constant-size steps was the only method which was able
to compute the deformed meshes.  It is also likely that UBN could compute 
the deformed meshes if the deformation were broken into smaller 
deformations (in a similar manner to the small-step FEMWARP algorithm 
described in~\cite{ShontzVava}).  We also showed
that UBN is much faster (i.e., up to 70 times faster) than the 
Newton continuation algorithm.  It could be argued that
continuation would be more competitive with UBN if only we had used
a different strategy for incrementing $\tau_k$.  This may be true,
but it seems to us that there is no good 
universal fast method for choosing the $\tau_k$. Our experiments
indicate, for example, that a more aggressive algorithm for
updating $\tau$ is more prone to terminating early due to
inverted elements.   Even selecting
the continuation path seems to be nontrivial (e.g., for the 2D annulus
example, it was necessary to parametrize the Dirichlet boundary condition
in polar rather than rectangular coordinates).  In contrast, the 
UBN method does not require any such analogous problem-dependent 
decisions, and the only parameters of the algorithm involve
termination criteria.

As described so far, our method applies to finite elements in which
the displacement field is piecewise linear over tetrahedra, but
UBN could be extended to piecewise quadratic
displacements.  The challenge with piecewise quad-ratic displacements
is that checking for tangling is much more complicated, as $J$ is not
constant on the element.  In particular, it is a function of both the
displacement and the location on the element, which makes it difficult
to determine when $J > 0$ analytically.   There are some
separate necessary and sufficient conditions for element inversion in the 
literature.  Let $G$ be the Jacobian of $F$.  
Then one such necessary condition 
is that $\det(G)$ has the same sign (strictly positive or strictly 
negative) 
at some finite list of test points~\cite{shephard}.  In this case, 
we would test for 
inversion at the Gauss points used for numerical quadrature; however, 
it is still 
possible that folding could occur at the corners.  A more complicated
sufficient condition 
for invertibility involving the Bernstein-B\'{e}zier form of a 
polynomial is given 
in~\cite{vavasis_bb}.  Checking that the sufficient condition is met 
requires running
a linear programming algorithm.  Salem, Canann, and Saigal have 
proposed sufficient 
conditions for quadratic triangles and tetrahedra 
in~\cite{salem1},~\cite{salem2},~\cite{salem3}, and~\cite{salem4}.

Another issue for UBN is uniqueness.  It should be noted that some classes 
of boundary value problems
may admit multiple solutions.  A somewhat complicated example of this
nonuniqueness
occurred in Section~\ref{sec:3Dexperiments}.  A conceptually
simpler example is as follows. Consider a long cylinder
in which one end is held at zero displacement, the other is rotated by
$2\pi$ radians, the long side-surface is traction-free, and
there are no body forces. Since the rotated nodes at one end return
to their original positions after rotation by $2\pi$ radians, a
valid solution to the boundary value problem is all zero displacements.
A second valid solution is a twisted configuration of the cylinder.

In the case of continuation, it is possible to select a sequence
of nonlinear boundary deformations to force the correct final configuration.
This is not possible with UBN, however, at least not without further
modification.
To distinguish one solution from another requires additional information
beyond boundary conditions.  Determining what form the additional
information ought to take will be studied as future work.
We will also determine 
how UBN should be extended in order to use such additional information.

Future work will also involve extending UBN to the incompressible or nearly
incompressible case.  In the incompressible case, the requirement that
$J=1$ becomes a constraint rather than a term in the energy
functional.  For this reason, $J$ disappears from the functional.  One
minimizes $\psi(\u)$ subject to the constraint $\g(\u)=\bz$, where the
latter expresses the $J=1$ constraint for each Gauss point.  The
functional $\psi(\u)$ typically involves the deviatoric strain at
Gauss points.  A common method for handling a constraint like this is
an augmented Lagrangian optimization algorithm \cite{NocedalWright}.
On each iteration of the augmented Lagrangian method, our UBN method is
applicable in the same way as in the unconstrained case considered
here.  In particular, the energy function $\psi(\u)$ is usually
undefined or nondifferentiable when $J=0$, so Newton's method is
unlikely to work well when $J$ gets close to zero or, even worse,
becomes negative.  Therefore, the preliminary untangling step and line
search described earlier are appropriate for the incompressible case
as well.

\section{Acknowledgements}
\label{sec:acknowledgements}
The authors wish to thank Patrick Knupp of Sandia National Laboratories 
for providing us with the 3D test meshes.  They benefited from many 
helpful conversations with Katerina Papoulia of University of Waterloo.
They also wish to thank the anonymous referees for their careful 
reading of the paper and for their helpful suggestions which strengthened 
it.

\end{document}